\documentclass[a4paper, 11pt]{article}
\usepackage{mathrsfs}

\usepackage{epsfig}
\usepackage{amsmath}
\usepackage{amssymb}
\usepackage{latexsym}
\usepackage{amsfonts}
\usepackage{amsthm}
\usepackage[numbers,sort&compress]{natbib}
\usepackage{graphicx}
\ifx\pdfoutput\undefined  \else
 \fi
\usepackage[centerlast]{caption2}
\usepackage{color}

\numberwithin{figure}{section} \numberwithin{definition}{section}
\numberwithin{observation}{section} \numberwithin{lemma}{section}
\numberwithin{theorem}{section} \numberwithin{proposition}{section}
\numberwithin{conjecture}{section}

\begin{document}

\title{
{\Large\bf Roman domination number of \\ Generalized Petersen Graphs
$P(n,2)$} \footnote{The research is supported by Chinese Natural
Science Foundations (60573022). }}
\author{
 Haoli Wang, \ Xirong Xu, \ Yuansheng Yang\footnote {corresponding
author's email : yangys@dlut.edu.cn} , Chunnian Ji \\
Department of Computer Science, \\Dalian University of
Technology\\
Dalian, 116024, P. R. China\\  }
\date{}
\maketitle
\begin{abstract}

\noindent A $Roman\ domination\ function$ on a graph $G=(V, E)$ is a
function $f:V(G)\rightarrow\{0,1,2\}$ satisfying the condition that
every vertex $u$ with $f(u)=0$ is adjacent to at least one vertex
$v$ with $f(v)=2$. The $weight$ of a Roman domination function $f$
is the value $f(V(G))=\sum_{u\in V(G)}f(u)$. The minimum weight of a
Roman dominating function on a graph $G$ is called the $Roman\
domination\ number$ of $G$, denoted by $\gamma_{R}(G)$. In this
paper, we study the {\it Roman domination number} of generalized
Petersen graphs $P(n,2)$ and prove that $\gamma_R(P(n,2)) = \lceil
{\frac{8n}{7}}\rceil (n \geq 5)$.


\end{abstract}

\section{Introduction}

\ \ \ \ Let $G=(V,E)$ be a simple graph, i.e., loopless and without
multiple edges, with vertex set $V(G)$ and edge set $E(G)$. The open
neighborhood, $N(v)$, and the closed neighborhood, $N[v]$, of a
vertex $v\in V$ are denoted by $N(v)=\{u \in V(G)\ :\ vu \in E(G)\}$
and $N[v]=N(v)\cup \{v\}$, respectively. For a vertex set $S
\subseteq V(G)$, $N(S)=\underset{v \in S}\cup N(v)$ and
$N[S]=\underset{v \in S}\cup N[v]$. The maximum degree of any vertex
in $V(G)$ is denoted by $\Delta(G)$.

A set $S \subseteq V(G)$ is a dominating set if for each $v \in
V(G)$ either $v \in S$ or $v$ is adjacent to some $w \in S$. That
is, $S$ is a dominating set if and only if $N[S]=V(G)$. The
domination number $\gamma(G)$ is the minimum cardinality of a
dominating set of $G$, and a dominating set $S$ of minimum
cardinality is called a $\gamma$-$set$ of $G$.

For a graph $G$, let $f:V\rightarrow\{0,1,2\}$, and let
$(V_0;V_1;V_2)$ be the order partition of $V$ induced by $f$, where
$V_i=\{v\in V(G)|f(v)=i\}$ and $|V_i|=n_i$, for $i=0,1,2$. Note that
there exists a 1-1 correspondence between the functions
$f:V(G)\rightarrow\{0,1,2\}$ and the ordered partitions
$(V_0;V_1;V_2)$ of $V(G)$. So we denote $f=(V_0;V_1;V_2)$.

A function $f:V(G)\rightarrow\{0,1,2\}$ is a $Roman\ domination\
function$ ($RDF$) if $V_2$ dominates $V_0$, i.e. $V_0\subseteq
N[V_2]$. The $weight$ of $f$ is $f(V(G))=\sum_{v\in
V(G)}f(v)=2n_2+n_1$. The minimum weight of an RDF of $G$ is called
the $Roman\ domination\ number$ of $G$, denoted by $\gamma_{R}(G)$.
We say that a function $f=(V_0;V_1;V_2)$ is a
$\gamma_{R}$-$function$ if it is an RDF and $f(V)=\gamma_{R}(G)$.

In 2004, Cockayne et al\cite{Cockayne04} studied the graph theoretic
properties of this variant of the domination number of a graph and
proved

\noindent{\bf Proposition 1.1. }For any graph $G$ of order $n$,
$\frac{2n}{\Delta(G)+1}\leq \gamma_R(G)$.

\noindent{\bf Proposition 1.2. }For any graph $G$ of order $n$,
$\gamma(G)\leq \gamma_R(G)\leq 2\gamma(G)$.

For more references and other Roman dominating problems, we can
refer to
\cite{Cockayne03,Hedetniemi03,Henning03,Henning02,ReVelle00,Stewart99}.

The generalized Petersen graph $P(n,k)$ is defined to be a graph on
$2n$ vertices with $V(P(n,k)) = \{v_i,u_{i}:0 \leq i \leq n-1\}$ and
$E(P(n,k)) = \{v_iv_{i+1},v_{i}u_{i},u_{i}u_{i+k}:0 \leq i \leq
n-1$, subscripts taken modulo $n \}$.

In 2007, Yang Yuansheng et al \cite{Yang06,Yang07} studied the
domination number of generalized Petersen graphs $P(n,2)$ and
$P(n,3)$. They proved

\noindent{\bf Theorem 1.3. } $\gamma$($P(n,2))=n - \lfloor
\frac{n}{5} \rfloor - \lfloor \frac{n+2}{5} \rfloor$.

\noindent{\bf Theorem 1.4. } $\gamma($P$(n,3))=
n-2\lfloor\frac{n}{4}\rfloor(n\neq 11)$.

In this paper, we study the $Roman$ domination in the generalized
Petersen graphs $P(n,2)$ and prove $\gamma_R$($P(n,2)) = \lceil
{\frac{8n}{7}} \rceil (n \geq 5)$.

\section{Roman domination number of
$P(n,2)$}

Let $m=\lfloor \frac{n}{7} \rfloor $, $t=n$(mod 7), then $n=7m+t$.

{\noindent \bf Lemma 2.1. }$\gamma_R$($P(n,2)) \leq
 \lceil {\frac{8n}{7}} \rceil (n \geq 5)$.

\noindent{\it Proof.} In order to prove that for $n \geq 5$,
$\gamma_R$($P(n,2)) \leq
 \lceil {\frac{8n}{7}} \rceil$, it suffices to give an RDF
$f$ of $P(n,2)$ with $f(V(P(n,2)))=\lceil {\frac{8n}{7}} \rceil$.
For $n=5$, let

$V_2=\{v_0,u_2,u_3\}$, $V_1=\emptyset$, $V_0=N(V_2)$.

\noindent For $n=6$, let

$V_2=\{v_0,u_3,u_4\}$, $V_1=\{v_2\}$, $V_0=N(V_2)$.

\noindent For $n\geq 7$, let

{\small
$$\begin{array}{llll}
V_2=\left \{\begin{array}{llll}
               \{v_{7i}, u_{7i+3}, u_{7i+4}\ :\ 0\leq i\leq m-1\}, & \mbox{if } t=0;\\
               \{v_{7i}, u_{7i+3}, u_{7i+4}\ :\ 0\leq i\leq m-1\}, & \mbox{if } t=1;\\
               \{v_{7i}, u_{7i+3}, u_{7i+4}\ :\ 0\leq i\leq m-1\}\cup \{v_{7m}\}, & \mbox{if } t=2;\\
               \{v_{7i}, u_{7i+3}, u_{7i+4}\ :\ 0\leq i\leq m-1\}\cup \{v_{7m}\}, & \mbox{if } t=3;\\
               \{v_{7i}, u_{7i+3}, u_{7i+4}\ :\ 0\leq i\leq
               m-1\}\cup \{v_{7m-1},u_{7m+1},u_{7m+2}\}, & \mbox{if } t=4;\\
               \{v_{7i}, u_{7i+3}, u_{7i+4}\ :\ 0\leq i\leq
               m-1\}\cup \{v_{7m},u_{7m+2},u_{7m+3}\}, & \mbox{if } t=5;\\
               \{v_{7i}, u_{7i+3}, u_{7i+4}\ :\ 0\leq i\leq m\}, & \mbox{if } t=6.
              \end{array}
           \right .\\
\\
V_1=\left \{\begin{array}{llll}
               \{v_{7i+2}, v_{7i+5}\ :\ 0\leq i\leq m-1\}, & \mbox{if } t=0;\\
               \{v_{7i+2}, v_{7i+5}\ :\ 0\leq i\leq m-1\} \cup \{v_{7m-1},u_{7m}\}, & \mbox{if } t=1;\\
               \{v_{7i+2}, v_{7i+5}\ :\ 0\leq i\leq m-1\} \cup \{u_{7m+1}\}, & \mbox{if } t=2;\\
               \{v_{7i+2}, v_{7i+5}\ :\ 0\leq i\leq m-1\} \cup \{u_{7m+1}, u_{7m+2}\}, & \mbox{if } t=3;\\
               \{v_{7i+2}, v_{7i+5}\ :\ 0\leq i\leq m-1\} \setminus \{v_{7m-2}\}, & \mbox{if } t=4;\\
               \{v_{7i+2}, v_{7i+5}\ :\ 0\leq i\leq m-1\}, & \mbox{if } t=5;\\
               \{v_{7i+2}, v_{7i+5}\ :\ 0\leq i\leq m-1\} \cup \{v_{7m+2}\}, & \mbox{if } t=6.
              \end{array}
           \right .\\
\\
V_0=N(V_2).
\end{array}$$}

Note that $V_0$, $V_1$ and $V_2$ are pairwise disjoint, and
$V(P(n,2))=V_1 \cup V_2 \cup V_0=V_1\cup N[V_2]$. Hence
$f=(V_0;V_1;V_2)$ is an RDF of $P(n,2)$ with {\small
$$\begin{array}{llll}
f(V(P(n,2)))=& \left \{\begin{array}{llll}
               2\times 3m+2m=8m=\lceil\frac{8n}{7}\rceil, & \mbox{ if } t=0;\\
               2\times 3m+2m+2=8m+2=\lceil\frac{8n}{7}\rceil, & \mbox{ if } t=1;\\
               2\times (3m+1)+2m+1=8m+3=\lceil\frac{8n}{7}\rceil, & \mbox{ if } t=2;\\
               2\times (3m+1)+2m+2=8m+4=\lceil\frac{8n}{7}\rceil, & \mbox{ if } t=3;\\
               2\times (3m+3)+2m-1=8m+5=\lceil\frac{8n}{7}\rceil, & \mbox{ if } t=4;\\
               2\times (3m+3)+2m=8m+6=\lceil\frac{8n}{7}\rceil, & \mbox{ if } t=5;\\
               2\times (3m+3)+2m+1=8m+7=\lceil\frac{8n}{7}\rceil, & \mbox{ if } t=6.
              \end{array}
           \right .
 \end{array}$$ }\hspace{200pt}\qed
\begin{figure}[ht]
\centering
\includegraphics[scale=0.9]{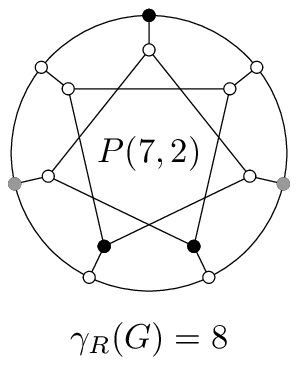}
\includegraphics[scale=0.9]{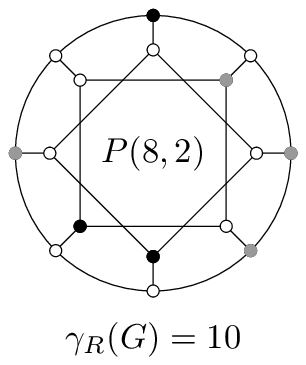}
\includegraphics[scale=0.9]{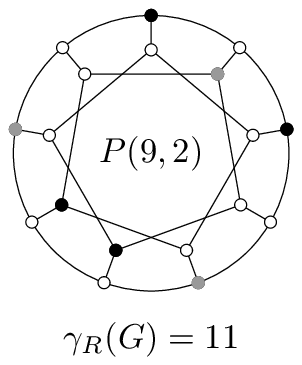}
\includegraphics[scale=0.9]{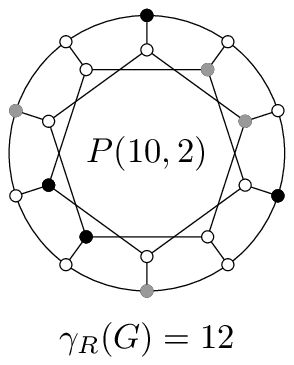}
\includegraphics[scale=0.9]{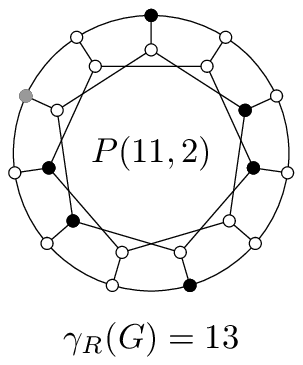}
\includegraphics[scale=0.9]{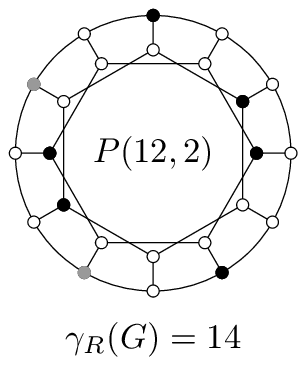}
\includegraphics[scale=0.9]{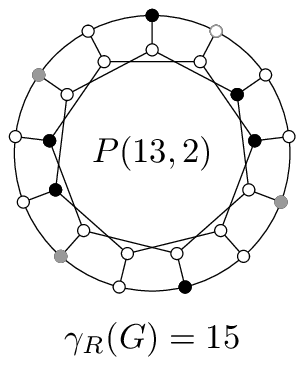}

\centering \vspace{5pt}\small{Figure 2.1:  The RDFs of $P(n,2)$ for
$7\leq n\leq 13$}
\end{figure}

In Figure 2.1, we give the RDFs of $P(n,2)$ for $7\leq n\leq 13$,
where the vertices of $V_2$ are in dark, the vertices of $V_1$ are
in grey, and the vertices of $V_0$ are in white.

Let $f=(V_0;V_1;V_2)$ be an arbitrary $\gamma_R$-function of
$P(n,2)$. Then we have

{\noindent \bf Lemma 2.2.}  For any vertex $w_1\in V_2$, if $w_2\in
N(w_1)$, then $w_2\not \in V_1$.

\noindent{\it Proof.} Suppose to the contrary that $w_2\in V_1$. Let
$f^{'}(w_2)=0$ and $f^{'}(w)=f(w)$ for every vertex $w\in
V(P(n,2))\setminus w_2$. Then $f^{'}$ is an RDF of $P(n,2)$ with
$f^{'}(V(P(n,2)))=\gamma_R(P(n,2))-1$, a contradiction.\qed

Let $f_m=(V_0;V_1;V_2)$ be an arbitrary $\gamma_R$-function of
$P(n,2)$ with minimum cardinality of $V_2$, i.e. $|V_2|\leq
|V^{'}_2|$ for any $\gamma_R$-function
$f^{'}=(V^{'}_0;V^{'}_1;V^{'}_2)$ of $P(n,2)$. Then we have

{\noindent \bf Lemma 2.3. } For any vertex $w_1\in V_2$, if $w_2\in
N(w_1)$, then $w_2\not \in V_2$.

\noindent{\it Proof.} Suppose to the contrary that $w_2\in
V(P(n,2))$ such that $w_1,w_2\in V_2$. Let $N(w_1)=\{w_2,w_3,w_4\}$
where $f_m(w_3)\geq f_m(w_4)$. There are two cases depending on
$w_4$:

\noindent{\bf Case 1.} $w_4\in V_1\cup V_2$. Let $f^{'}(w_1)=0$ and
$f^{'}(w)=f_m(w)$ for every vertex $w\in V(P(n,2))\setminus
\{w_1\}$. Then $f$ is an RDF of $P(n,2)$ with
$f^{'}(V(P(n,2)))=\gamma_R(P(n,2))-2$, a contradiction.

\noindent{\bf Case 2.} $w_4\in V_0$. For $w_3\in V_1\cup V_2$, let
$f^{'}(w_1)=0$, $f^{'}(w_4)=1$ and $f^{'}(w)=f_m(w)$ for every
vertex $w\in V(P(n,2))\setminus \{w_1,w_4\}$. We have that $f$ is an
RDF of $P(n,2)$ with $f^{'}(V(P(n,2)))=\gamma_R(P(n,2))-1$, a
contradiction. For $w_3\in V_0$, let $f^{'}(w_1)=0$,
$f^{'}(w_3)=f^{'}(w_4)=1$ and $f^{'}(w)=f_m(w)$ for every vertex
$w\in V(P(n,2))\setminus \{w_1,w_3,w_4\}$. We have that
$f'=(V^{'}_0;V^{'}_1;V^{'}_2)$ is a $\gamma_R$-function of $P(n,2)$
with $|V_2|>|V^{'}_2|$, a contradiction. \qed

For an arbitrary $\gamma_R$-function $f=(V_0;V_1;V_2)$ of $P(n,2)$,
we define a function $g_f$ as follows. Let
$$\begin{array}{llll}
g_{f}(w) = & \left \{\begin{array}{llll}
               0.5, & $ if $ w\in V_2;\\
               1, & $ if $ w\in V_1;\\
               0.5|N(w)\cap V_2|, & $ if $ w\in V_0.
               \end{array}
           \right .
\end{array}$$
Then $g_{f}(w)\geq 0.5$ for every vertex $w\in V(P(n,2))$.

{\noindent \bf Lemma 2.4.} $g_{f}(V(P(n,2)))=\sum_{v\in
V(P(n,2))}g_f(v)=\gamma_R(P(n,2))$.

\noindent{\it Proof.} By Lemmas 2.2-2.3, we have that $N(w)\subseteq
V_0$ for any vertex $w\in V_2$. It follows that
$\gamma_R(P(n,2))=|V_1|+2|V_2|=|V_1|+0.5|V_2|+0.5\sum_{w\in
V_2}|N(w)\bigcap V_0|=g_{f}(V_1)+g_{f}(V_2)+0.5\sum_{w\in
V_0}|N(w)\bigcap V_2|=g_{f}(V_1)$
$+g_{f}(V_2)+g_f(V_0)=g_{f}(V(P(n,2)))$. \qed

For every vertex $w\in V(P(n,2))$, let $r_f(w)=g_f(w)-0.5$. For
every subset $S\subseteq V(P(n,2))$, let $r_{f}(S)=\sum_{w\in
S}(r_f(w))$. Let $V^{'}(i,t)=\{v_j,u_j:i\leq j\leq i+t-1\}$. Then we
have

{\noindent \bf Lemma 2.5.} If $r_f(V^{'}(i,7))\leq 0.5$, then
$v_{i+3}\not \in V_2$.

\noindent{\it Proof.} Suppose to the contrary that $v_{i+3}\in V_2$.
Then $v_{i+2},v_{i+4},u_{i+3}\in V_0$. If $\{v_{i+5},u_{i+5}\}\cap
(V_1\cup V_2)\neq \emptyset$, then $rf(V^{'}(i,7))=0.5$. It follows
that $v_{i+1},u_{i+1}\in V_0$, $v_i\in V_2$, $u_i,u_{i+2}\in V_0$
and $u_{i+4}\in V_2$, which implies that $r_f(V^{'}(i,7))\geq 1$(see
Figure 2.2(1) for $v_{i+5}\in V_1$), a contradiction. Hence,
$v_{i+5},u_{i+5}\in V_0$. It follows that $v_{i+6}\in V_2$. There
are three cases depending on $u_{i+4}$:

\noindent \textbf{Case 1.} $u_{i+4}\in V_2$. Then
$r_f(V^{'}(i,7))\geq 1$(see Figure 2.2(2)), a contradiction.

\noindent \textbf{Case 2.} $u_{i+4}\in V_1$. Then
$r_f(V^{'}(i,7))\geq 0.5$. It follows $u_{i+2}\in V_0$, $u_i\in V_2$
and at least one vertex of $\{v_{i+1},u_{i+1}\}$ belongs to $V_1\cup
V_2$, which implies that $r_f(V^{'}(i,7))\geq 1$(see Figure 2.2(3)),
a contradiction.

\noindent \textbf{Case 3.} $u_{i+4}\in V_0$. Then $u_{i+2}\in V_2$
and $r_f(V^{'}(i,7))\geq 0.5$. Since at least one vertex of
$\{v_i,v_{i+1},u_{i+1}\}$ belongs to $V_1\cup V_2$, we have
$r_f(V^{'}(i,7))\geq 1$(see Figure 2.2(4)), a contradiction.

From the above discussion, the lemma follows. \qed

\begin{figure}[ht]
\centering
\includegraphics[scale=0.6]{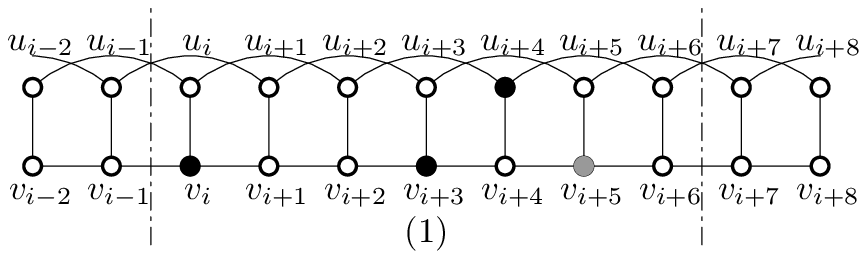}\hspace{20pt}
\includegraphics[scale=0.6]{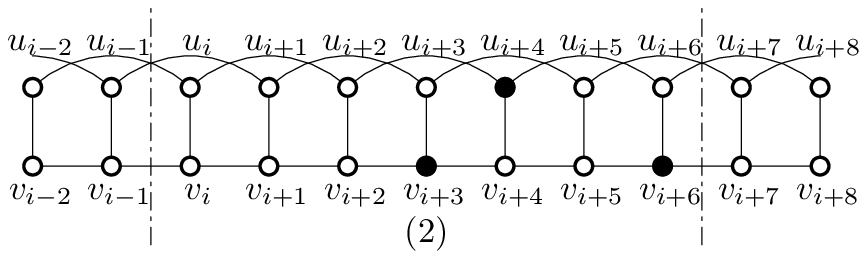}
\includegraphics[scale=0.6]{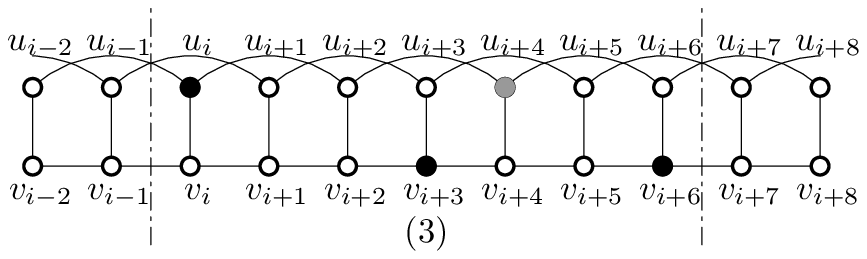}\hspace{20pt}
\includegraphics[scale=0.6]{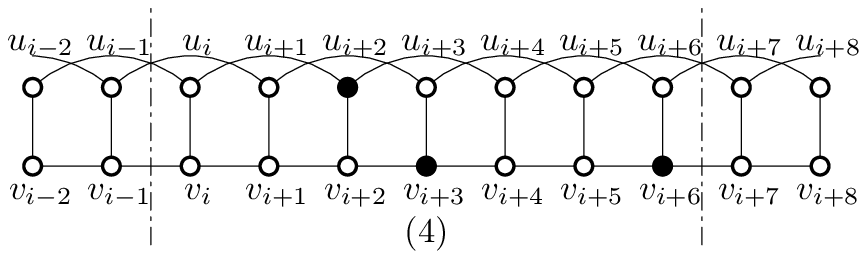}

\centering\small{Figure 2.2:  Some cases for $v_{i+3}\in V_2$ }
\end{figure}

{\noindent \bf Lemma 2.6.} If $r_f(V^{'}(i,7))\leq 0.5$, then
$v_{i+2},v_{i+4}\not \in V_2$.

\noindent{\it Proof.} By symmetry, it suffices to prove that
$v_{i+4}\not \in V_2$. Suppose to the contrary that $v_{i+4}\in
V_2$. Then $v_{i+3},v_{i+5},u_{i+4}\in V_0$. If
$\{u_{i+5},v_{i+6},u_{i+6}\}\cap (V_1\cup V_2)\neq \emptyset$, then
$rf(V^{'}(i,7))= 0.5$. It follows that $v_{i+2},u_{i+2}\in V_0$,
$v_{i+1}\in V_2$, $u_{i+2}\in V_0$ and $u_{i}\in V_2$, which implies
that $r_f(V^{'}(i,7))\geq 1$(see Figure 2.3(1) for $u_{i+5}\in
V_1$), a contradiction. Hence, $u_{i+5},v_{i+6},u_{i+6}\in V_0$.
Then $v_{i+7}\in V_2$, $u_{i+3}\in V_2$ and $r_f(V^{'}(i,7))\geq
0.5$. It forces $u_{i+1},v_{i+1},v_{i+2}\in V_0$, $v_{i}\in V_2$ and
$N[u_{i+2}]\subseteq V_0$(see Figure 2.3(2)), a contradiction. \qed

\begin{figure}[ht]
\centering
\includegraphics[scale=0.6]{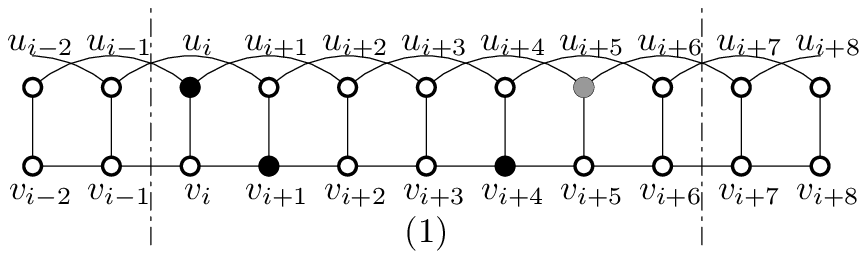}\hspace{20pt}
\includegraphics[scale=0.6]{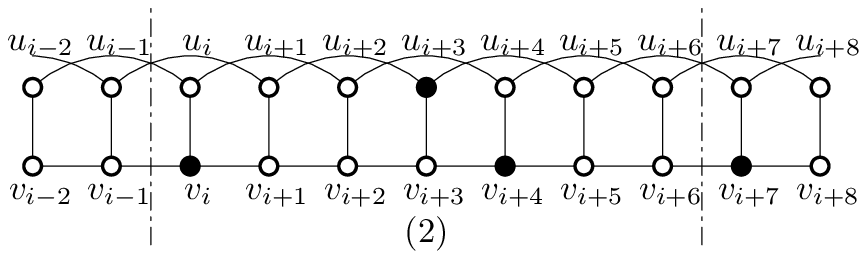}

\centering\small{Figure 2.3:  Some cases for $v_{i+4}\in V_2$ }
\end{figure}

{\noindent \bf Lemma 2.7.} If $r_f(V^{'}(i,7))\leq 0.5$, then
$u_{i+3}\not \in V_2$.

\noindent{\it Proof.} Suppose to the contrary that $u_{i+3}\in V_2$.
Then $u_{i+1},v_{i+3},u_{i+5}\in V_0$. If $\{v_{i+4},v_{i+5}\}\cap
(V_1\cup V_2)\neq \emptyset$, then $r_f(V^{'}(i+4,3))= 0.5$. If
$v_{i+4},v_{i+5}\in V_0$, then $u_{i+4},v_{i+6}\in V_2$, we also
have $r_f(V^{'}(i+4,3))= 0.5$. It follows that $v_{i+1}\in V_0$ and
$v_{i+2}\in V_1$, which implies that $r_f(V^{'}(i+4,3))\geq 1$(see
Figure 2.4), a contradiction. \qed

\begin{figure}[ht]
\centering
\includegraphics[scale=0.6]{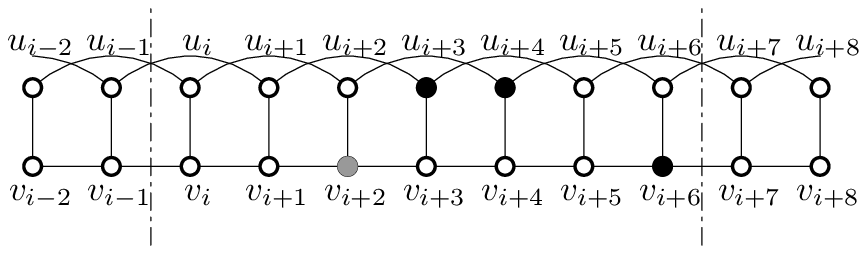}

\centering\small{Figure 2.4:  Case for $u_{i+3}\in V_2$ }
\end{figure}

{\noindent \bf Lemma 2.8.} If $r_f(V^{'}(i,7))\leq 0.5$, then \\
\indent (1) $r_f(V^{'}(i,7))$ is at least 0.5,\\
\indent (2) $V^{'}(i,7)\cap V_1=\{v_{i+3}\}$,\\
\indent (3) $V^{'}(i,7)\cap V_2=\{u_{i+1},u_{i+2},v_{i+5}\}$ and
$V^{'}(i,7)\cap V_0=V^{'}(i,7)\setminus
\{v_{i+3},u_{i+1},u_{i+2},v_{i+5}\}$, or $V^{'}(i,7)\cap
V_2=\{u_{i+4},u_{i+5},v_{i+1}\}$ and $V^{'}(i,7)\cap
V_0=V^{'}(i,7)\setminus \{v_{i+3},u_{i+4},u_{i+5},v_{i+1}\}$.

\noindent{\it Proof.} By Lemma 2.5-2.7, we have that
$r_f(V^{'}(i,7))$ is at least 0.5 and $v_{i+3}\in V_1$. It follows
that $u_{i+3},v_{i+2},v_{i+4}\in V_0$ and one vertex of
$\{u_{i+1},u_{i+5}\}$ belongs to $V_2$(see Figure 2.5(1)). If
$u_{i+1}\in V_2$, then $u_{i+2},v_{i+5}\in V_2$. If $u_{i+5}\in
V_2$, then $u_{i+4},v_{i+1}\in V_2$(see Figure 2.5(2)). \qed

\begin{figure}[ht]
\centering
\includegraphics[scale=0.6]{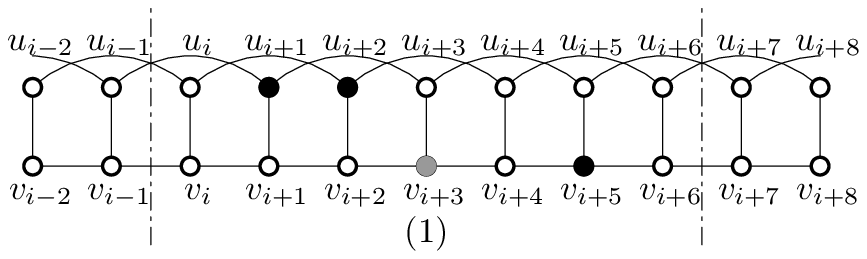} \hspace{20pt}
\includegraphics[scale=0.6]{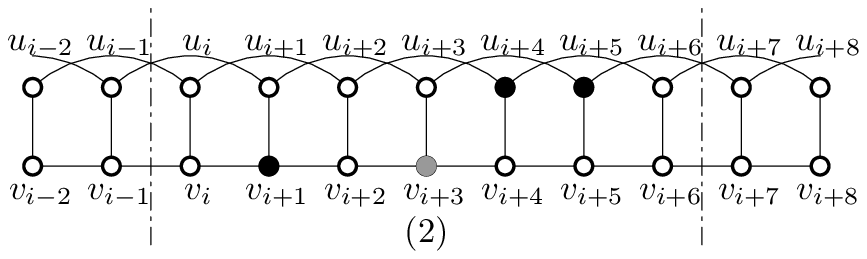}

\centering\small{Figure 2.5:  Case for $v_{i+3}\in V_1$ }
\end{figure}

{\noindent \bf Lemma 2.9.} If $r_f(V^{'}(i,7))= 0.5$, then
$r_f(V^{'}(i-7,7))\geq 1$ and $r_f(V^{'}(i+7,7))\geq 1.5$, or
$r_f(V^{'}(i-7,7))\geq 1.5$ and $r_f(V^{'}(i+7,7))\geq 1$.

\noindent{\it Proof.} By symmetry, we only need to consider the case
shown in Figure 2.5(1). Since $N(v_i)\cap V_2\neq \emptyset$ and
$N(u_{i+6})\cap V_2\neq \emptyset$, we have $v_{i-1},u_{i+8}\in
V_2$. By Lemma 2.5 and Lemma 2.7, we have $r_f(V^{'}(i-4,7))\geq 1$
and $r_f(V^{'}(i+5,7))\geq 1$. We see that $r_f(V^{'}(i,3))=0$ and
$r_f(V^{'}(i+5,2))=0$. It follows that
$r_f(V^{'}(i-4,4))=r_f(V^{'}(i-4,7))-r_f(V^{'}(i,3))\geq 1$ and
$r_f(V^{'}(i+7,5))=r_f(V^{'}(i+5,7))-r_f(V^{'}(i+5,2))\geq 1$.
Therefore, $r_f(V^{'}(i-7,7))\geq 1$ and $r_f(V^{'}(i+7,7))\geq 1$.

Now, we prove that $r_f(V^{'}(i-7,7))\geq 1.5$. Suppose to the
contrary that $r_f(V^{'}(i-7,7))=1$. Since $r_f(V^{'}(i,7))=0.5$, we
have $u_{i-2}\not \in V_2$. If $u_{i-2}\in V_1$, then
$r_f(V^{'}(i-7,7))\geq 1$ and $u_{i-3},v_{i-3}\in V_0$. It follows
that $u_{i-5},v_{i-4}\in V_2$ and $r_f(V^{'}(i-7,7))\geq 1.5$(see
Figure 2.5(1)), a contradiction. Hence, $u_{i-2}\in V_0$. It follows
that $u_{i-4}\in V_2$ and $v_{i-4}\in V_0$. There are three cases
depending on $u_{i-3}$:

\noindent{\bf Case 1.} $u_{i-3}\in V_2$. Then $r_f(V^{'}(i-7,7))\geq
1$ and $u_{i-5}\in V_0$. It follows that at least one vertex of
$\{v_{i-6},v_{i-5}\}$ belongs to $V_1\cup V_2$, which implies that
$r_f(V^{'}(i-7,7))\geq 1.5$(see Figure 2.5(2)), a contradiction.

\noindent{\bf Case 2.} $u_{i-3}\in V_1$. Then $v_{i-3}\in V_1$,
which implies that $r_f(V^{'}(i-7,7))\geq 1.5$(see Figure 2.5(3)), a
contradiction.

\noindent{\bf Case 3.} $u_{i-3}\in V_0$. Then $v_{i-3}\in V_1$ and
$u_{i-5}\in V_2$. It follows that at least one vertex of
$\{v_{i-7},v_{i-6}\}$ belongs to $V_1\cup V_2$, which implies that
$r_f(V^{'}(i-7,7))\geq 1.5$(see Figure 2.5(4)), a contradiction.

From the above discussion, the lemma follows. \qed
\begin{figure}[ht]
\centering
\includegraphics[scale=0.6]{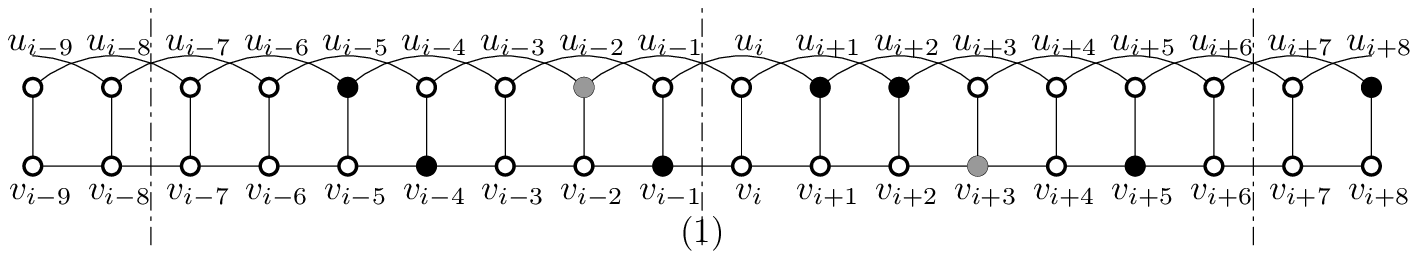}
\includegraphics[scale=0.6]{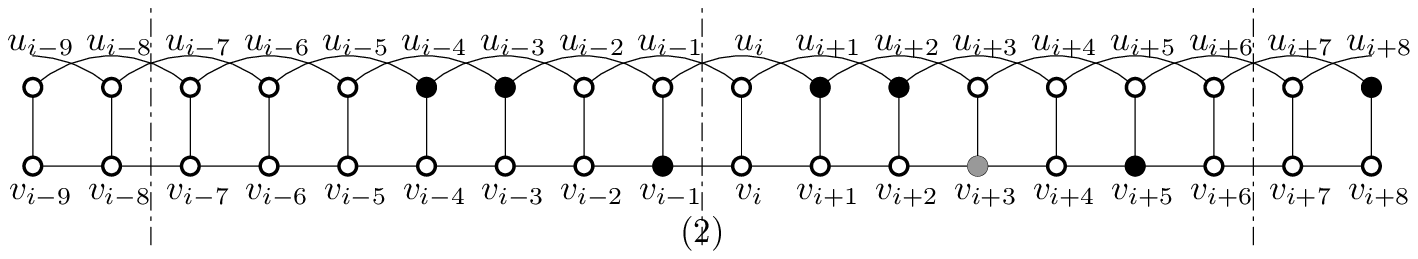}
\includegraphics[scale=0.6]{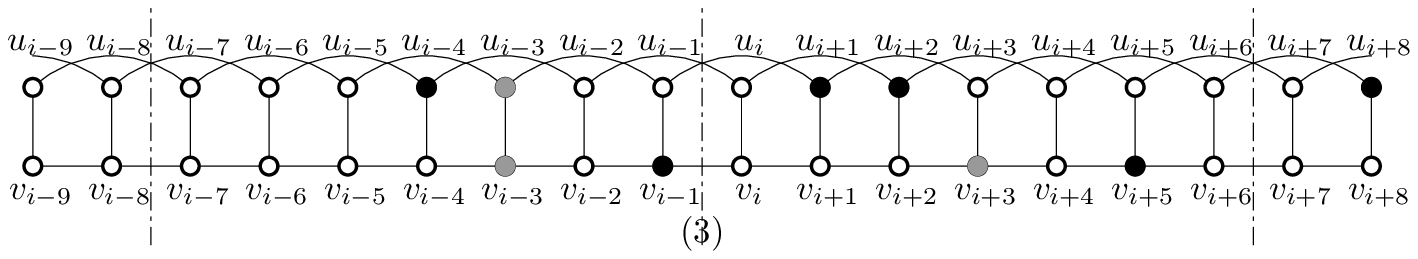}
\includegraphics[scale=0.6]{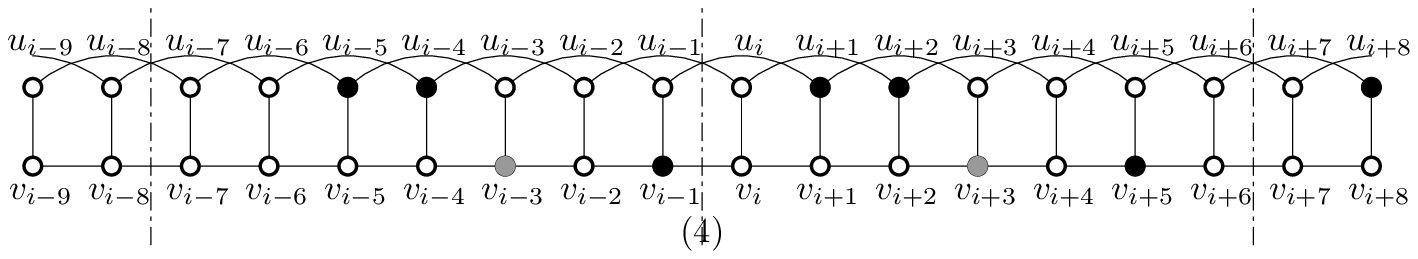}

\centering\small{Figure 2.5:  Some cases for $r_f(V^{'}(i-7,7))=1$ }
\end{figure}

{\noindent \bf Lemma 2.10.} If $r_f(V^{'}(i-7,7))= 0.5$ and
$r_f(V^{'}(i+7,7))= 0.5$, then $r_f(V^{'}(i,7))\geq 2$.

\noindent{\it Proof.} By Lemma 2.8, there are three cases.

\noindent{\bf Case 1.} $v_{i-2},u_{i+8}\in V_2$. Then
$u_{i+1},v_{i+6}\in V_2$. It follows that $v_{i+1},u_{i+3}$,
$v_{i+5},u_{i+6}\in V_0$ and $r_f(V^{'}(i,7))\geq 0.5$. Since
$N[u_{i+5}]\cap (V_1\cup V_2)\neq \emptyset$, we have $u_{i+5}\in
V_1\cup V_2$, which implies that $r_f(V^{'}(i,7))\geq 1$. Since
$r_f(V^{'}(i-7,7))= 0.5$, we have $v_i,u_i\not\in V_2$. Since
$N[v_i]\cap (V_1\cup V_2)\neq \emptyset$, we have $v_i\in V_1$. Then
$r_f(V^{'}(i,7))\geq 1.5$. Since $N[v_{i+3}]\cap (V_1\cup V_2)\neq
\emptyset$, we have $\{v_{i+2},v_{i+3},v_{i+4}\}\cap (V_1\cup
V_2)\neq \emptyset$, which implies that $r_f(V^{'}(i,7))\geq 2$ (see
Figure 2.6).

\begin{figure}[ht]
\centering
\includegraphics[scale=0.57]{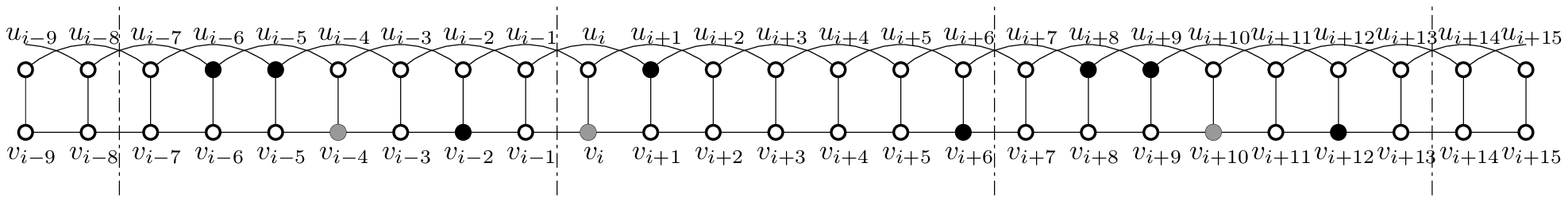}

\centering\small{Figure 2.6:  Case for $v_{i-2},u_{i+8}\in V_2$ }
\end{figure}

\noindent{\bf Case 2.} $v_{i-2},v_{i+8}\in V_2$. Then
$u_{i+1},u_{i+5}\in V_2$. It follows that $v_{i+1},u_{i+3}$,
$v_{i+5}\in V_0$ and $r_f(V^{'}(i,7))\geq 0.5$. Since $N[v_i]\cap
(V_1\cup V_2)\neq \emptyset$, we have $v_i\in V_1$. Since
$N[v_{i+6}]\cap (V_1\cup V_2)\neq \emptyset$, we have $v_{i+6}\in
V_1$. It follows that $r_f(V^{'}(i,7))\geq 1.5$. Since
$N[v_{i+3}]\cap (V_1\cup V_2)\neq \emptyset$, we have
$\{v_{i+2},v_{i+3},v_{i+4}\}\cap (V_1\cup V_2)\neq \emptyset$, which
implies $r_f(V^{'}(i,7))\geq 2$ (see Figure 2.7).

\begin{figure}[ht]
\centering
\includegraphics[scale=0.57]{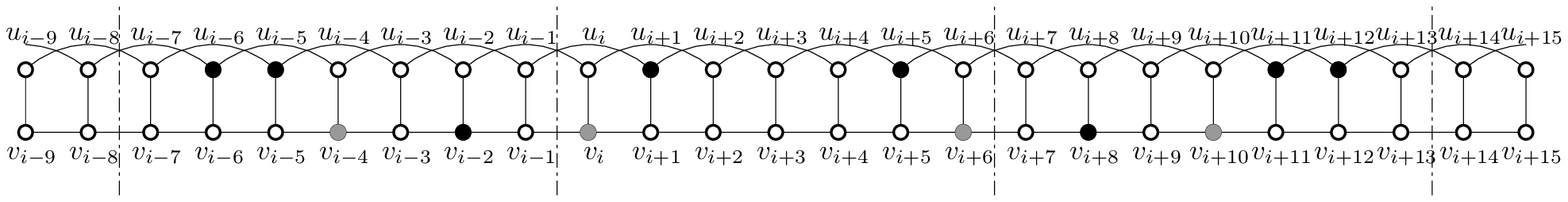}

\centering\small{Figure 2.7:  Case for $v_{i-2},v_{i+8}\in V_2$ }
\end{figure}

\noindent{\bf Case 3.} $u_{i-2},u_{i+8}\in V_2$. Then
$v_i,v_{i+6}\in V_2$. It follows that $u_i,v_{i+1},v_{i+5}$,
$u_{i+6}\in V_0$ and $r_f(V^{'}(i,7))\geq 1$. Since
$r_f(V^{'}(i-7,7))= 0.5$ and $r_f(V^{'}(i+7,7))= 0.5$, we have
$u_{i+1},u_{i+5}\not\in V_2$.

If $u_{i+1}\in V_1$ and $u_{i+5}\in V_1$, then
$r_f(V^{'}(i+7,7))\geq 2$ (see Figure 2.8(1)).

If $u_{i+1}\in V_0$ or $u_{i+5}\in V_0$, then $u_{i+3}\in V_2$. It
follows that $u_{i+1},v_{i+3}$, $u_{i+5}\in V_0$. Since
$N[v_{i+2}]\cap (V_1\cup V_2)\neq \emptyset$ and $N[v_{i+4}]\cap
(V_1\cup V_2)\neq \emptyset$, we have that $\{v_{i+2},u_{i+2}\}\cap
(V_1\cup V_2)\neq \emptyset$ and $\{v_{i+4},u_{i+4}\}\cap  (V_1\cup
V_2)\neq \emptyset$, which implies that $r_f(V^{'}(i+7,7))\geq 2$
(see Figure 2.8(2)).

This completes the proof. \qed

\begin{figure}[ht]
\centering
\includegraphics[scale=0.57]{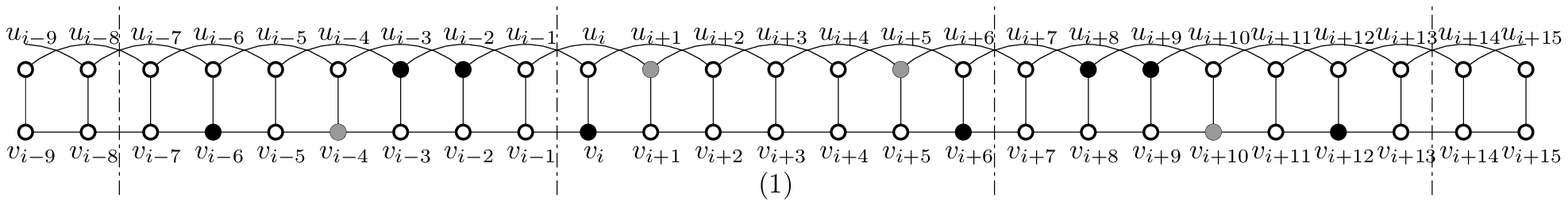}
\includegraphics[scale=0.57]{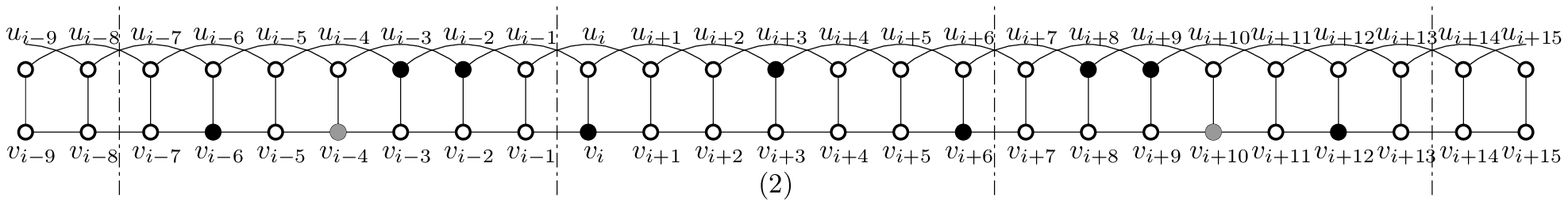}

\centering\small{Figure 2.8:  Case for $u_{i-2},u_{i+8}\in V_2$ }
\end{figure}

{\noindent \bf Lemma 2.11.} $\gamma_R(P(n,2)) \geq
 \lceil {\frac{8n}{7}} \rceil (n \geq 5)$.

\noindent{\it Proof.} Let
$$\begin{array}{llll}
S_1&=\{i:0\leq i\leq n-1,r_f(V'(7i,7))=0.5\},\\
S_2&=\{i:0\leq i\leq n-1,r_f(V'(7i,7))=1\},\\
 S_{31}&=\{i:0\leq i\leq n-1,r_f(V'(7i,7))\geq
1.5$, $|\{i-1,i+1\}\cap S_1|\leq 1\},\\
S_{32}&=\{i:0\leq i\leq n-1,r_f(V'(7i,7))\geq 1.5,\ |\{i-1,
i+1\}\cap S_1|=2\}.
\end{array}$$
By Lemma 2.8, $r_f(V'(7i,7))\geq 0.5$, hence we have
$\{0,1,\ldots,n-1\}=S_1\cup S_2\cup S_{31}\cup S_{32}$. By Lemma
2.9, we have $|S_1|\leq |S_{31}|+2|S_{32}|$. By Lemma 2.10, we have
that $r_f(V'(7i,7))\geq 2$ for any integer $i\in S_{32}$. By Lemma
2.4, we have $$\begin{array}{llll}
&7\times \gamma_R(P(n,2))\\
=&7\times \sum_{v\in V(P(n,2))}g_f(v)\\
=&7\times \sum_{v\in V(P(n,2))}(r_f(v)+0.5)\\
=&7\times \sum_{v\in V(P(n,2))}r_f(v)+7n\\
=&\sum_{0\leq i\leq n-1}r_f(V'(7i,7))+7n\\
=&\sum_{i\in S_1}r_f(V'(7i,7))+\sum_{i\in
S_2}r_f(V'(7i,7))+\sum_{i\in
S_{31}}r_f(V'(7i,7))\\&+\sum_{i\in S_{32}}r_f(V'(i,7))+7n\\
\geq&0.5|S_1|+|S_2|+1.5|S_{31}|+2|S_{32}|+7n\\
=&0.5|S_1|+|S_2|+|S_{31}|+|S_{32}|+0.5(|S_{31}|+2|S_{32}|)+7n\\
\geq& |S_1|+|S_2|+|S_{31}|+|S_{32}|+7n\\
=&8n,
\end{array}$$ which implies that $\gamma_R(P(n,2))\geq \lceil \frac{8n}{7}\rceil$.
\hspace{200pt}\qed

From Lemma 2.1 and Lemma 2.11, we have the following

{\noindent \bf Theorem 2.12. }$\gamma_R(P(n,2)) = \lceil
{\frac{8n}{7}}\rceil (n \geq 5)$. \hspace{14pt}$\Box{}$

\end{document}